\documentclass[final]{siamltex}

\usepackage{amsfonts,amssymb}
\usepackage{amssymb}
\usepackage{cmmib57}

\newtheorem{remark}[theorem]{{\it Remark}}
\newtheorem{example}[theorem]{{\it Example}}

\def\bf{\bfseries}
\def\it{\itshape}

\newcommand{\eps}{\varepsilon}

\newcommand{\text}[1]{\hbox{\rm \ #1\ \/}}

\newcommand{\bl}[1]{\begin{lemma}\label{#1}}
\newcommand{\br}[1]{\begin{remark}\label{#1}}
\newcommand{\bt}[1]{\begin{theorem}\label{#1}}
\newcommand{\bd}[1]{\begin{definition}\label{#1}}
\newcommand{\bp}[1]{\begin{proposition}\label{#1}}
\newcommand{\bc}[1]{\begin{corollary}\label{#1}}
\newcommand{\bfact}[1]{\begin{fact}\label{#1}}
\newcommand{\bex}[1]{\begin{example}\label{#1}}
\newcommand{\bem}[1]{\begin{example}\label{#1}}
\newcommand{\ec}{\end{corollary}}
\newcommand{\eex}{\end{example}}
\newcommand{\eem}{\end{example}}
\newcommand{\el}{\end{lemma}}
\newcommand{\er}{\end{remark}}
\newcommand{\et}{\end{theorem}}
\newcommand{\ed}{\end{definition}}
\newcommand{\ep}{\end{proposition}}
\newcommand{\epr}{\end{proof}}
\newcommand{\bpr}{\begin{proof}}

\newcommand{\beq}{\begin{eqnarray}}
\newcommand{\eeq}{\end{eqnarray}}
\newcommand{\beqn}{\begin{eqnarray*}}
\newcommand{\eeqn}{\end{eqnarray*}}
\newcommand{\bi}{\begin{itemize}}
\newcommand{\ei}{\end{itemize}}
\newcommand{\ben}{\begin{enumerate}}
\newcommand{\een}{\end{enumerate}}

\newcommand{\R}{{\mathbb R}}  
\newcommand{\N}{{\mathbb N}}  

\newcommand{\ds}{\displaystyle}

\newcommand{\cK}{{\cal K}}
\newcommand{\cKL}{{\cal KL}}
\newcommand{\sgn}{{\rm sgn}}

\newcommand{\calo}{{\cal O}}

\title{Global Asymptotic Controllability Implies Input to State
Stabilization\thanks{Received by the editors...  }}
\author{Michael Malisoff\thanks{Corresponding Author.  Department of
Mathematics, 304 Lockett Hall, Louisiana State University and A \&
M College,
 Baton Rouge LA 70803-4918, USA
(malisoff@lsu.edu).  This author was supported in part by
Louisiana Board of Regents Support Fund Contract
LEQSF(2002-04)-ENH-TR-26, as part of the project
``Interdisciplinary Education, Outreach, and Research in Control
Theory at LSU''.} \and Ludovic Rifford\thanks{Institut Girard
Desargues, Universit\'e Lyon 1, B\^atiment Braconnier, 21 Avenue
Claude Bernard, 69622 Villeurbanne Cedex, France
(rifford@desargues.univ-lyon1.fr).} \and Eduardo
Sontag\thanks{Department of Mathematics, Rutgers-New Brunswick,
Hill Center-Busch Campus, 110 Frelinghuysen Road, Piscataway NJ
08854-8019, USA (sontag@control.rutgers.edu).   This author was
supported
 by USAF Grant F49620-01-1-0063, and by NSF Grant
CCR-0206789.}}
\begin{document}

\maketitle \vspace{-1.36in} \slugger{sicon}{200?}{??}{?}{???--???}
\vspace{1.06in}

\setcounter{page}{1}
\begin{abstract}
The main problem addressed in this paper is the design of
feedbacks for globally asymptotically controllable (GAC) control
affine systems that render the closed loop systems input to state
stable with respect to actuator errors. Extensions for fully
nonlinear GAC systems with actuator errors are also discussed. Our
controllers have the property that they tolerate small observation
noise as well.
\end{abstract}



\begin{PII}

\end{PII}

\pagestyle{myheadings} \thispagestyle{plain} \markboth {MICHAEL
MALISOFF, LUDOVIC RIFFORD, EDUARDO SONTAG} {CONTROLLABILITY AND
INPUT TO STATE STABILIZATION}

\section{Introduction}
\label{introd}

The theory of input to state stability (ISS) forms the basis for
much of current research in mathematical control theory (see
\cite{PW96,S99, S00}).  The ISS property was introduced in
\cite{S89a}.  In the past decade, there has been a great deal of
research done on the problem of finding ISS stabilizing control
laws (see  \cite{FK96, KKK95, KL98, LSW01}).
 This note is
concerned with the ISS of control systems of the form
\begin{equation}\label{dyn}
\dot x=f(x)+G(x)u
\end{equation}
where $f$ and $G$ are locally Lipschitz vector fields on $\R^n$,
$f(0)=0$, and the control $u$ is valued in $\R^m$ (but see also
$\S$\ref{rks} for extensions for fully nonlinear systems). We
assume throughout that (\ref{dyn}) is globally asymptotically
controllable (GAC), and we construct a feedback $K:\R^n\to\R^m$
for which
\begin{equation}
\label{twoprime} \dot x=f(x)+G(x)K(x)+G(x)u
\end{equation}
is ISS.  As pointed out in \cite{B83, SS80}, a continuous
stabilizing feedback $K$ fails to exist in general.  This fact
forces us to consider discontinuous feedbacks $K$, so our
solutions will be taken in the more general sense of sampling and
Euler solutions for dynamics that are discontinuous in the state.
By an Euler solution, we mean a uniform limit of sampling
solutions, taken as the frequency of sampling becomes infinite
(see $\S$\ref{defs} for precise definitions). This will extend
\cite{S89a, S90}, which show how to make $C^o$-stabilizable
systems ISS to actuator errors.  In particular, our results apply
to the nonholonomic integrator  (see  \cite{B83, LR03}, and
$\S$\ref{fbrk} below) and other applications where Brockett's
condition is not satisfied, and which therefore cannot be
stabilized by continuous feedbacks (see  \cite{S98, S99, St02}).

Our results  also {\em strengthen} \cite{CLSS97}, which
constructed feedbacks for GAC systems that render the closed loop
systems globally asymptotically stable.  Our main tool will be the
recent constructions of semiconcave control Lyapunov functions
(CLF's) for GAC systems from \cite{R00, R02}. Our results  also
apply in the more general situation where measurement noise may
occur. In particular, our feedback $K$ will have the additional
feature that the {\em perturbed} system
\begin{equation}
\label{two}
 \dot x= f(x)+G(x) K(x+e)+G(x) u
\end{equation}
is also ISS when the observation error $e:[0,\infty)\to \R^n$ in
the controller
 is
{\em sufficiently small}.  In this context, the precise value of
$e(t)$ is unknown to the controller, but information about upper
bounds on the magnitude of $e(t)$ can be used to design the
feedback. We will prove the following:

\bt{thm1} If (\ref{dyn}) is GAC, then there exists a feedback $K$
for which  (\ref{two}) is ISS for Euler solutions.\et

The preceding theorem characterizes the uniform limits of sampling
solutions of (\ref{two}) (see  $\S$\ref{defs} for the precise
definitions of Euler and sampling solutions). From a computational
standpoint, it is also desirable to know how frequently to sample
in order to achieve ISS for sampling solutions.  This information
is provided in the following semi-discrete version of Theorem
\ref{thm1} for sampling solutions:

\bt{thm2} If (\ref{dyn}) is GAC, then there exists a feedback $K$
for which  (\ref{two}) is ISS for sampling solutions.\et

This paper is organized as follows.  In $\S$\ref{defs}, we review
the relevant background on   CLF's, ISS, nonsmooth analysis, and
discontinuous feedbacks. In $\S$\ref{proofs}, we prove our main
results. This is followed in $\S$\ref{fbrk} by a comparison of our
feedback construction with the known feedback constructions for
$C^o$-stabilizable systems, and an application of our results to
the nonholonomic integrator. We close in $\S$\ref{rks} with an
extension for fully nonlinear systems.

\section{Definitions and Main Lemmas}
\label{defs}
 Let $\cK_\infty$ denote the set of all continuous functions
$\rho:[0,\infty)\to[0,\infty)$ for which (i) $\rho(0)=0$ and  (ii)
$\rho$ is strictly increasing and unbounded.  Note for future
reference that $\cK_\infty$ is closed under inverse and
composition (i.e., if $\rho_1,\rho_2\in \cK_\infty$, then
$\rho^{-1}_1, \rho_1\circ\rho_2\in \cK_\infty$). We  let $\cKL$
denote the set of all continuous functions $\beta:[0,\infty)\times
[0,\infty)\to[0,\infty)$ for which (1) $\beta(\cdot, t)\in {\cal
K}_\infty$ for each $t\ge 0$, (2) $\beta(s,\cdot)$ is
nonincreasing for each $s\ge 0$, and (3) $\beta(s,t)\to 0$ as
$t\to +\infty$ for each $s\ge 0$.

For each $k\in \N$ and $r>0$, we define \[{\cal M}^k=\{{\rm
measurable\ } u:[0,\infty)\to \R^k: |u|_\infty <\infty\}\] and
${\cal M}^k_r:=\{u\in {\cal M}^k: |u|_\infty\le r\}$, where
$|\cdot|_\infty$ is the essential supremum.  We let $\Vert
u(s)\Vert_I$ denote the essential supremum of a function $u$
restricted to an interval $I$.  Let $|\cdot |$ denote the
Euclidean norm, in the appropriate dimension, and
\[r{\cal B}_k:=\{x\in\R^k: |x|< r\}\] for each $k\in \N$ and $r>0$.  The
closure of $r{\cal B}_k$ is denoted by $r\bar {\cal B}_k$, and
${\rm bd} (S)$ denotes the boundary of any subset $S$ in Euclidean
space. We also set
\[{\cal O}:=\{e:[0,\infty)\to\R^n\}, \; \; \sup(e)=\sup\{|e(t)|:
t\ge 0\}\] for all $e\in {\cal O}$, and ${\cal O}_\eta:=\{e\in
{\cal O}: \sup(e)\le \eta\}$ for each $\eta>0$. For any compact
set ${\cal F}\subseteq \R^n$ and $\eps>0$, we define the compact
set
\[{\cal F}^\eps:=\{x\in\R^n: \min\{|x-p|: p\in {\cal F}\}\le \eps\},\]
i.e., the ``$\eps$-enlargement of $\cal F$''.   Given a continuous
function \[h:\R^n\times \R^m\to \R^n: (x,u)\mapsto h(x,u)\] that
is locally Lipschitz in $x$ uniformly on compact subsets of
$\R^n\times \R^m$, we let $\phi_h(\cdot, x_o,u)$ denote the
trajectory of $\dot x=h(x,u)$ starting at $x_o\in\R^n$ for each
choice of  $u\in {\cal M}^{m}$. In this case, $\phi_h(\cdot,
x_o,u)$ is defined on some maximal interval $[0, t)$, with $t>0$
depending on $u$ and $x_o$.  Let $C^k$ denote the set of all
continuous functions $\varphi:\R^n\to\R$ that have at least $k$
continuous derivatives (for $k=0,1$).
 We use the following controllability notion, which was introduced in
\cite{S83} and later reformulated in terms of ${\cal KL}$
functions in \cite{S99}:

\begin{definition}\label{gacdef} We call the system $\dot x=h(x,u)$  {\rm
globally asymptotically controllable (GAC)} provided  there are a
nondecreasing function  $\sigma:[0,\infty)\to[0,\infty)$ and a
function $\beta\in {\cal KL}$ satisfying the following: For each
$x_o\in\R^n$, there exists $u\in {\cal M}^m$ such that \bi
\item[(a)] $|\phi_h(t,x_o,u)|\le \beta(|x_o|,t)$ for all $t\ge 0$;
\ and \item[(b)] $|u(t)|\le \sigma(|x_o|)$ for a.e. $t\ge 0$. \ei
In this case, we call $\sigma$ the {\rm GAC modulus} of $\dot
x=h(x,u)$. \ed

In our main results, the controllers will be taken to be
discontinuous feedbacks, so the dynamics will be discontinuous in
the state variable. Therefore, we will form our trajectories
through sampling,  and through uniform limits of sampling
trajectories, as follows. We say that
$\pi=\{t_o,t_1,t_2,\ldots\}\subset [0,\infty)$ is a {\it
partition} of $[0,\infty)$ provided $t_o=0$, $t_i < t_{i+1}$ for
all $i\ge 0$, and $t_i\to +\infty$ as $i\to +\infty$. The set of
all partitions of $[0,\infty)$ is denoted by ${\rm Par}$. Let
\[F:\R^n\times \R^m\times \R^m\to \R^n: (x,p,u)\mapsto F(x,p,u)\] be
a continuous function that is locally Lipschitz in $x$ uniformly
on compact subsets of $\R^n\times \R^m\times \R^m$.  A {\it
feedback} for $F$ is defined to be any locally bounded function
$K:\R^n\to\R^m$ for which $K(0)=0$.  In particular, we allow
discontinuous feedbacks. The arguments $x$, $p$, and $u$ in $F$
are used to represent the state, feedback value, and actuator
error, respectively.

  Given a feedback
$K:\R^n\to\R^m$, $\pi=\{t_o, t_1,t_2,\ldots\}\in {\rm Par}$,
$x_o\in \R^n$, $e\in {\cal O}$, and $u\in {\cal M}^m$, the {\it
sampling solution} for the initial value problem (IVP)
\begin{eqnarray}
\dot x(t) &=& F(x(t),K(x(t)+e(t)),u(t))\label{feed}\\
x(0)&=&x_o\label{side}
\end{eqnarray}
 is the continuous function defined by recursively solving
\begin{equation}
\label{recu} \dot x(t) \; \; =\; \;
F(x(t),K(x(t_i)+e(t_i)),u(t))\end{equation} from the initial time
$t_i$ up to time $s_i=t_i\vee\sup\{s\in [t_i,t_{i+1}]: x(\cdot)
{\rm\ is\ defined\ on\ } [t_i,s)\}$, where $x(0)=x_o$. In this
case, the sampling solution of (\ref{feed})-(\ref{side}) is
defined on the right-open interval from time zero up to time $\bar
t=\inf\{s_i: s_i<t_{i+1}\}$.
    This sampling solution will be
denoted by $t\mapsto x_\pi(t; x_o,u,e)$ to exhibit its dependence
on $\pi\in {\rm Par}$, $x_o\in\R^n$, $u\in {\cal M}^m$, and $e\in
\calo$, or simply by $x_\pi$, when the dependence is clear from
the context.  Note that if $s_i=t_{i+1}$ for all $i$, then $\bar
t=+\infty$ (as the infimum of the empty set), so in that case, the
sampling solution $t\mapsto x_\pi(t; x_o,u,e)$ is defined on
$[0,\infty)$.

 We also define the
{\it  upper diameter} and the {\it lower diameter} of a given
partition $\pi=\{t_o,t_1,t_2,\ldots\}$ by
\[\ds\overline{\mathbf d}(\pi)=\sup_{i\ge 0}(t_{i+1}-t_i),\; \; \;
\ds\underline{\mathbf d}(\pi)=\inf_{i\ge 0}(t_{i+1}-t_i)\]
respectively.   We let ${\rm Par}(\delta):=\left\{\pi\in {\rm
Par}: \overline{\mathbf d}(\pi)<\delta\right\}$ for each
$\delta>0$.
 We will say that a function $y: [0,\infty)\to \R^n$ is an
 {\it  Euler solution (robust to small observation errors)}
of
\begin{equation}
\label{store} \dot x(t)=F(x(t),K(x(t)), u(t)),\; \; x(0)=x_o
\end{equation}
 for $u\in {\cal M}^m$ provided there are
sequences $\pi_r\in {\rm Par}$ and  $e_r\in {\cal O}$ such that
\bi \item[(a)] $\overline{\mathbf d}(\pi_r)\to 0$; \item[(b)]
$\sup(e_r)/\underline{\mathbf d}(\pi_r)\to 0$; \ and \item[(c)]
$t\mapsto x_{\pi_r}(t; x_o, u, e_r)$ converges uniformly to $y$ as
$r\to +\infty$. \ei Note that the approximating trajectories in
the preceding definition  all use the same input $u$ (but see
Remark~\ref{euleremark} for a more general notion of Euler
solutions, which also involves sequences of inputs).

This paper will design feedbacks  that make closed loop GAC
systems ISS with respect to actuator errors. More precisely, we
will use the following definition:

\bd{Defn3} We
 say that  (\ref{feed}) is {\rm ISS for sampling
 solutions}
provided there are $\beta\in {\cal KL}$ and $\gamma\in {\cal
K}_\infty$ satisfying:  For each $\eps, M,N>0$ with $0<\eps<M$,
there exist positive $\delta=\delta (\eps, M,N)$ and
$\kappa=\kappa(\eps, M,N)$ such that for each $\pi\in {\rm
Par}(\delta)$, $x_o\in M\bar {\cal B}_n$, $u\in {\cal M}^m_N$, and
$e\in {\cal O}$ for which $\sup(e)\le \kappa \underline{\mathbf
d}(\pi)$,
\begin{equation}\label{siss}
|x_\pi(t; x_o, u,e)|\le \max\{\beta(M,t)+\gamma(N),\eps\}
\end{equation}
 for
all $t\ge 0$.\ed

Roughly speaking, condition (\ref{siss}) says that the system is
ISS,
     modulo small overflows, if the sampling is done `quickly enough', as
     determined by the condition $\pi\in {\rm Par}(\delta)$, but `not too
     quickly', as determined by the additional requirement that $\underline{\mathbf
     d}(\pi)\ge(1/\kappa)\sup(e)$.  In the special case where the observation
     error $e\equiv 0$, the condition on $\underline{\mathbf d}(\pi)$ in
     Definition~\ref{Defn3} is no longer needed; our results are new even for
     this particular case.

  Notice that the bounds on $e$
are in the supremum, not the essential supremum. It is easy to
check that Definition~\ref{Defn3} remains unchanged if we replace
the right-hand side in (\ref{siss}) by
$\beta(M,t)+\gamma(N)+\eps$. We also use the following analog of
Definition~\ref{Defn3} for Euler solutions:

\bd{Defn2} We say that the system (\ref{feed}) is {\rm ISS for
Euler solutions} provided there exist $\beta\in {\cal KL}$ and
$\gamma\in {\cal K}_\infty$ satisfying:  If $u\in {\cal M}^m$ and
$x_o\in\R^n$, and if $t\mapsto x(t)$ is an Euler solution of
(\ref{store}), then
\begin{equation}\label{iss}
|x(t)|\le
\beta(|x_o|,t)+\gamma(|u|_\infty)\end{equation}
for all $t\ge 0$.
\end{definition}

\begin{remark}\label{euleremark}\rm In the definition of Euler
solutions we gave above, all of the approximating trajectories
$t\mapsto x_{\pi_r}(t; x_o, u, e_r)$ use the same input $u\in
{\cal M}^m$. A different way to define Euler solutions, which
gives rise to a more general class of limiting solutions, is as
follows: A function $y: [0,\infty)\to \R^n$ is a
 {\it  generalized Euler solution}
of (\ref{store}) for $u\in {\cal M}^m$ provided there are
sequences $\pi_r\in {\rm Par}$, $e_r\in {\cal O}$, and $u_r\in
{\cal M}^m$ such that \bi \item[(a)] $\overline{\mathbf
d}(\pi_r)\to 0$; \item[(b)] $\sup(e_r)/\underline{\mathbf
d}(\pi_r)\to 0$; \item[(c)] $|u_r|_\infty\le |u|_\infty$ for all
$r$;\  and \item[(d)] $t\mapsto x_{\pi_r}(t; x_o, u_r, e_r)$
converges uniformly to $y$ as $r\to +\infty$. \ei We can then
define ISS for generalized Euler solutions exactly as in
Definition~\ref{Defn2}, by merely replacing ``Euler solution''
with ``generalized Euler solution'' throughout the definition. Our
proof of Theorem~\ref{thm1} will actually show the following
slightly more general result:  If (\ref{dyn}) is GAC, then there
exists a feedback $K$ for which (\ref{two}) is ISS for generalized
Euler solutions. \er

 Our main tools in this paper will be
nonsmooth analysis and nonsmooth Lyapunov functions.  The
following definitions will be used.  Let $\Omega$ be an arbitrary
open subset of $\R^n$.  Recall the following definition:

\bd{scc} Let $g:\Omega\to \R$ be a continuous function on
$\Omega$;  it is said to be {\rm semiconcave} on $\Omega$ provided
for each point $x_o\in \Omega$, there exist $\rho, C>0$ such that
\[
g(x)+g(y)-2g\left(\frac{x+y}{2}\right)\le C||x-y||^2
\]
for all $x,y\in x_o+\rho {\cal B}_n$. \ed

The {\it proximal superdifferential} (respectively, {\it  proximal
subdifferential}) of a function $V: \Omega\to \R$ at $x\in
\Omega$, which is denoted by $\partial^PV(x)$ (resp., $\partial_P
V(x)$), is defined to be the set of all $\zeta\in\R^n$ for which
there exist $\sigma,\eta>0$ such that
\[V(y)-V(x)-\sigma|y-x|^2\le \langle \zeta, y-x\rangle\; \; \;
{\rm (resp.,}\; V(y)-V(x)-\sigma|y-x|^2\ge \langle \zeta,
y-x\rangle)\] for all $y\in x+\eta {\cal B}_n$.  The {\it limiting
subdifferential} of a continuous function $V:\Omega\to\R$ at
$x\in\Omega$ is
\[\partial_LV(x):=\{q\in\R^n:   \exists x_n\to x \; \; \&\; \;  q_n\in
\partial
_PV(x_n) \; \; {\rm s.t.}\; \;  q_n\to q\}.\]

In what follows, we assume  $h:\R^n\times \R^m\to \R^n:
(x,u)\mapsto h(x,u)$ is continuous, that it is locally Lipschitz
in $x$ uniformly on compact subsets of $\R^n\times \R^m$, and that
$h(0,0)=0$. The following definition was introduced in \cite{S83}
and reformulated in proximal terms in \cite{S99}:
 \bd{Defn5} A {\rm
control-Lyapunov function (CLF)} for
\begin{equation}\label{feq}
\dot x=h(x,u)
\end{equation}
is a continuous, positive definite, proper function $V:\R^n\to \R$
for which there exist a continuous, positive definite function
$W:\R^n\to \R$, and a nondecreasing function
$\alpha:[0,\infty)\to[0,\infty)$,  satisfying
\[
\ds \forall \zeta\in \partial_PV(x),\; \; \inf_{|u|\le
\alpha(|x|)}\langle \zeta , h(x,u)\rangle\le -W(x)
\]
for all $x\in\R^n$.  In this case, we call $(V,W)$ a {\rm Lyapunov
pair} for {\rm (\ref{feq})}. \ed

Recall the following lemmas (see  \cite{R02}):

\begin{lemma}\label{lemma2} If (\ref{feq}) is GAC, then there exists a CLF
$V$ for (\ref{feq}) that is semiconcave on $\R^n\setminus\{0\}$,
and a nondecreasing function $\alpha:[0,\infty)\to [0,\infty)$,
that  satisfy
\begin{eqnarray}
\label{eleven} \forall \zeta \in \partial_L V(x), \; \;
\min_{|u|\le \alpha(|x|)} \langle \zeta, h(x,u)\rangle \leq -V(x)
\end{eqnarray}
for all $x\in\R^n$.
\end{lemma}

\begin{lemma}\label{lemma3}  Let $V:\Omega\to\R$ be semiconcave. Then $V$
is locally Lipschitz, and
 $\emptyset\ne \partial_LV(x)\subseteq
\partial^PV(x)$ for all $x\in\Omega$.  Moreover, for each compact
 set $Q\subset \Omega$, there exist constants $\sigma, \mu>0$ such that
$V(y)-V(x)-\sigma |y-x|^2\le \langle \zeta,y-x\rangle $ for all
$y\in x+\mu {\cal B}_n$, all $x\in Q$, and all $\zeta\in
\partial^PV(x)$.

\end{lemma}

Notice that Lemma~\ref{lemma3} allows the constants  in the
definition of $\partial^PV(x)$ to be chosen uniformly on compact
sets.

\begin{remark}\rm In \cite{R02}, the controls $u$ take all their values in
a given compact metric space $U$. The precise version of the CLF
existence theorem in  \cite{R02} is the same as our Lemma
\ref{lemma2}, except that the infimum in the decay condition
(\ref{eleven}) is replaced by the infimum over all $u\in U$.
  The version of
Lemma~\ref{lemma2} we gave above follows from a slight
modification of the arguments of \cite{R00, R02}, using the GAC
modulus in the GAC definition (see Definition~\ref{gacdef}). The
existence theory \cite{R00} for semiconcave CLF's is a
strengthening of the proof that continuous CLF's exist for any GAC
system (see  \cite{S83}).
 \er

\section{Proofs of Theorems}
\label{proofs}

Let $V$ be a CLF satisfying the requirements of  Lemma
\ref{lemma2} for the dynamics
\begin{equation}
\label{sis}
 h(x,u)=f(x)+G(x)u.
\end{equation}
Define the  functions $\underline{\alpha},\overline{\alpha}\in
{\cal K}_\infty$ by
\begin{equation}
\label{checkcheck} \underline{\alpha}(s)=\min\{|x|: V(x)\ge s\}\;
\; {\rm and}\; \;  \overline{\alpha}(s)=\max\{|x|: V(x)\le s\}.
\end{equation}
One can easily check that
\begin{equation}
\label{easily} \forall x\in\R^n,\; \; \;
\underline{\alpha}(V(x))\le |x|\; \; {\rm and}\; \;
\overline{\alpha}(V(x))\ge |x|.
\end{equation}
Moreover, by reducing $\underline \alpha$, we may assume that
$\underline\alpha(s)\le s$ for all $s\ge 0$, while still
satisfying (\ref{easily}).

 Let $x\mapsto \zeta(x)$
be any selection of $\partial _LV(x)$ on $\R^n\setminus\{0\}$,
with $\zeta(0):=0$.   For each $x\in\R^n$, we can   choose
$u=u_x\in \alpha(|x|){\cal B}_m$ that satisfies the inequality in
(\ref{eleven}) for the dynamics (\ref{sis}) and $\zeta=\zeta(x)$.
Define the feedback $K_1:\R^n\to\R^m$ by $K_1(x)=u_x$ for all
$x\ne 0$ and $K_1(0)=0$. We use the functions
\begin{equation}
\label{newww}
\begin{array}{l}
a(x)=\langle \zeta(x),f(x)+G(x)K_1(x)\rangle,\; \;
b_j(x)= \langle \zeta(x),g_j(x)\rangle \; \forall j\\
K_2(x) =-V(x) (\sgn\{b_1(x)\}, \sgn\{b_2(x)\},\ldots,
\sgn\{b_m(x)\})^T \end{array},\end{equation}
 where $g_j$ is the $j$th column of
$G$ for $j=1,2,\ldots, m$, and
\[
\sgn\{s\}\; \; =\; \;  \left\{
\begin{array}{rl}
1,& s>0\\
-1,& s<0\\
0,& s=0
\end{array}
\right..\] We remark that our results remain true, with minor
changes in the proofs, if the factor $-V(x)$ in the definition of
$K_2$ is replaced by $-W(x)$ for an arbitrary positive definite
proper continuous function $W:\R^n\to \R$.
In particular, $K:=K_1+K_2$ is a feedback for the dynamics
\[
F(x,p,u)=f(x)+G(x)(p+u).
\]
Moreover,
\begin{equation}\label{aV}
a(x)\le -V(x)<0\; \; \forall x\in \R^n\setminus\{0\}.
\end{equation}
We next show that
\begin{equation}
\label{feed2} \dot x(t)= F(x(t), K(x(t)+e(t)),u(t))
\end{equation}
is ISS for sampling solutions.

To this end, choose $\eps, M,N>0$ for which $0<\eps<M$.  It
clearly suffices to verify the ISS property (\ref{siss}) for
$\eps<1$, since that would imply the property for {\em all}
overflows $\eps>0$. Choose
\begin{equation}\label{digs}
u\in {\cal M}^m_N,\; \;  e\in {\cal O}_{\eps/16},\; \;  x_o\in
M\bar {\cal B}_n.
\end{equation}
 In what follows, $x_\pi$ denotes the sampling solution
for (\ref{feed2}) for the choices (\ref{digs}) and $\pi\in {\rm
Par}$, and $\tilde x_\pi$ is the (possibly discontinuous) function
that is inductively defined by solving the IVP \[ \dot
x(t)=f(x(t))+G(x(t))[K(\tilde x_i)+u(t)],\; \; x(t_i)=\tilde x_i
\] on $[t_i,t_{i+1})$, where $\tilde x_i
:=x_i+e(t_i)$, $x_i:=x_\pi(t_i)$,  and $\pi=\{t_o,
t_1,t_2,\ldots\}$.  We later restrict the choice of $\pi$ so that
$x_\pi$ and $\tilde x_\pi$ are defined on $[0,\infty)$.  We will
use the compact set \[ Q=\left\{\left[\overline{\alpha}\circ
\underline{\alpha}^{-1}(N+M)+1\right]\bar {\cal
B}_n\right\}\setminus \eps {\cal B}_n.\]
 Notice that $Q, Q^{\eps/2}\subseteq \R^n \setminus\{0\}$, and
 that $x_o\in Q^\eps$.
 Using
Lemma~\ref{lemma3} and the semiconcavity of $V$ on
$\R^n\setminus\{0\}$, we can find $\sigma, \mu>0$ such that
\begin{equation}
\label{oo} V(y)-V(x)\le \langle \zeta(x), y-x\rangle+\sigma
|y-x|^2
\end{equation}
for all $y\in x+\mu {\cal B}_n$ and $x\in Q^{\eps/2}$. Let ${\cal
L}_\eps>1$ be a Lipschitz constant for $V$ on $Q^{\eps/2}$, the
existence of which is also guaranteed by Lemma~\ref{lemma3}. It
follows from the definition of a CLF that
\begin{equation}\label{check}
\begin{array}{l}
\lambda_-:=\min\left\{V(p): p\in Q^{\eps/2}\right\}\\
\lambda_+:=\max\left\{V(p): p\in Q^{\eps}\right\}
\end{array}
\end{equation}
are finite positive numbers.  Therefore, we can choose
$\tilde\eps\in (0,\eps)$ for which
\begin{equation}\label{stars}
\overline{\alpha}\left(p+\frac{{\cal
L}_\eps}{4}\tilde\eps\right)\le
\overline{\alpha}(p)+\frac{\eps}{8}\; \; \; \; \; \; \forall p\in
\left[0,\underline{\alpha}^{-1}(N)+\lambda_+\right].
\end{equation}
We can also  find
\begin{equation}
\label{circlecross} \delta=\delta(\eps, M,N)\in
\left(0,\frac{\tilde\eps}{16+\lambda_++16\lambda_+}\right)\end{equation}
 such that
if
\begin{equation}
\label{bigsmall} \pi\in {\rm Par}(\delta),\; \;  e\in {\cal
O}_{\tilde\eps/16}, \; \; x_i\in Q^\eps,\end{equation} and if
$t\in [t_i,t_{i+1})$ is such that $x_\pi(s)$ and $\tilde x_\pi(s)$
remain in $Q^{2\eps}$ for all $s\in [t_i,t]$, then
\begin{equation}
\label{shrp} \max\{|x_\pi(t)-x_i|,
 |\tilde x_\pi(t)-\tilde x_i|\}\le
 \min\left\{\mu,\frac{\tilde\eps}{16(1+{\cal L}_\eps)},
\sqrt{\frac{\lambda_-}{8\sigma}(t-t_i)}\right\}.
\end{equation}
 This follows from the local boundedness of
$K$, $f$ and $G$. It follows from (\ref{shrp}) that $\tilde
x_\pi(t)\in Q^{\eps/4}$ (resp., $x_\pi(t)\in Q^{\eps/4}$) for all
$t\in [t_i,t_{i+1})$ and all $i$ such that $\tilde x_i\in Q$
(resp., $x_i\in Q$), since the trajectories cannot move the
initial value more than $\frac{\eps}{16}$ and there are no blow up
times for the trajectories. In particular, (\ref{shrp}) will show
that $x_\pi$ and $\tilde x_\pi$ are defined on $[0,\infty)$, since
the argument we are about to give shows that $x_i\in Q^\eps$ for
all $i$. By reducing $\delta$ as necessary, we can assume
\begin{eqnarray}\label{thund}
 \ds
 \begin{array}{l}
\Vert\zeta(\tilde x_i)\cdot \left(F(\tilde x_i, K(\tilde x_i),
u(s))- f(\tilde x_\pi(s))\right.\\\left.-G(\tilde
x_\pi(s))[u(s)+K(\tilde x_i)]\right)\Vert_{[t_i, t_{i+1})}\le
\frac{\lambda_-}{8}\end{array}
\end{eqnarray}
 for all $i$ such that $\tilde x_i\in
Q^{\eps/2}$. This follows from  the  Lipschitzness of $f$ and $G$
on $Q^\eps$. Having chosen $\delta$ to satisfy the preceding
requirements, pick any $\pi\in {\rm Par}(\delta)$. It follows from
(\ref{oo}) and (\ref{shrp}) that
\begin{equation}
\label{diff}
\begin{array}{ccl}
 V(\tilde x_\pi(t))-V(\tilde x_i)&\le& \langle
\zeta(\tilde x_i), \tilde x_\pi(t)-\tilde x_i\rangle+\sigma
|\tilde x_\pi(t)-\tilde x_i|^2\\ &\le& \langle\zeta(\tilde x_i),
\tilde x_\pi(t)-\tilde x_i\rangle + \frac{\lambda_-}{8}(t-t_i)
\end{array}\end{equation}
for all $t\in [t_i,t_{i+1})$ and all $i$ such that $\tilde x_i\in
Q^{\eps/4}$. Moreover, if $\tilde x_i\in Q^{\eps/4}$ and $t\in
[t_i,t_{i+1})$, and if
\begin{equation}\label{xcirc}
V(\tilde x_i)\ge N,
\end{equation}
 then
\begin{eqnarray}
\label{bigg}\langle \zeta(\tilde x_i),\tilde x_\pi(t)-\tilde
x_i\rangle &\le & \left\langle \zeta(\tilde x_i),\int_{t_i}^t
F(\tilde x_i, K(\tilde x_i), u(s)) {\rm
d}s\right\rangle+\frac{\lambda_-}{8}(t-t_i)\; \;\; \; ({\rm by}\;
(\ref{thund})) \nonumber\\\nonumber &=& (t-t_i) \langle
\zeta(\tilde x_i),f(\tilde x_i)+G(\tilde x_i)K(\tilde x_i) \rangle
\\\nonumber
&+& \int_{t_i}^t \langle \zeta(\tilde x_i),G(\tilde x_i)
u(s)\rangle {\rm d}s +\frac{\lambda_-}{8}(t-t_i)
\\\label{rou1}
&\le & (t-t_i) a(\tilde x_i)-(t-t_i)V(\tilde
x_i)\sum_{j=1}^m|b_j(\tilde x_i)|\\\nonumber &+& N(t-t_i)
\sum_{j=1}^m|b_j(\tilde x_i)| +\frac{\lambda_-}{8}(t-t_i)
\\\nonumber
&\le & (t-t_i)a(\tilde x_i)+\frac{\lambda_-}{8}(t-t_i) \; \; \;
({\rm by}\; (\ref{xcirc}))
\\\nonumber
&\le & -(t-t_i)V(\tilde x_i)+\frac{\lambda_-}{8}(t-t_i)\; \; \;
({\rm by}\; (\ref{aV})).
\end{eqnarray}
 Let
\[S=\{x\in\R^n: V(x)\le \underline{\alpha}^{-1}(N)\}.\]
Then $S\subset Q^\eps$.  Indeed, $x\in S$ implies
\[
\underline{\alpha}\circ \overline{\alpha}^{-1}(|x|)\le
\underline{\alpha}\circ \overline{\alpha}^{-1}\circ
\overline{\alpha}\circ V(x)\le N,\] and therefore  $|x|\le
\overline{\alpha}\circ\underline{\alpha}^{-1}(N)$.  By further
reducing $\eps$, we can assume  $(2\eps){\cal B}_n\subset S$. If
$\tilde x_i\in Q^{\eps/4}$ but $\tilde x_i\not\in S$, then
$V(\tilde x_i)\ge \underline\alpha^{-1}(N)\ge N$, so (\ref{check})
and (\ref{diff}) give
\begin{eqnarray}
\label{circc}
 V(\tilde x_\pi(t))-V(\tilde x_i)&\le &
-(t-t_i)\frac{V(\tilde x_i)}{2}+(t-t_i)\frac{\lambda_-}{4}\nonumber\\
 &\le & -(t-t_i)\frac{V(\tilde x_i)}{4}\; \; \forall
 t\in[t_i,t_{i+1}).
\end{eqnarray}
Let ${\cal L}_f$ and ${\cal L}_G$  be  Lipschitz constants for $f$
and $G$ restricted to  $Q^\eps$, respectively. Define the
constants
\begin{equation}
\label{constants}
\begin{array}{l}
R=N+\sup\left\{|K(x)|: x\in Q^{\eps/2}\right\},\\ L={\cal
L}_f+R{\cal L}_G,\; \; \;
\ds\kappa=\kappa(\eps,M,N):=\frac{\min\{\lambda_-,\eps\}}{16{\cal
L}_\eps(e^{L\delta}+1)}.
\end{array}
\end{equation}
We will presently show that
\begin{equation}\label{circstar}
\ds \sup_{t_i\le t<t_{i+1}} |x_\pi(t)-\tilde x_\pi(t)|\le
|e(t_i)|e^{L\delta}\; \; \; \; \; \forall i \; \; {\rm s.t.}\; \;
x_i\in Q^{\eps/4}.
\end{equation}
Using (\ref{circstar}), we will now find $\beta\in{\cal KL}$ and
$\gamma\in \cK$ to satisfy the ISS estimate
\begin{equation}
\label{goal} |x_\pi(t)|\le \beta(|x_o|,t)+\gamma(N)+\eps\; \;
\forall t\ge 0
\end{equation}
which will prove Theorem~\ref{thm2}.

To this end, assume $x_i\in Q$, but that $x_i\not\in
S^{\tilde\eps/16}$. Then (\ref{shrp}) implies $x_\pi(t)$ and
$\tilde x_\pi(t)$ both remain in $Q^{\eps/4}$ on $[t_i,t_{i+1})$.
Moreover, $\tilde x_i\in Q^{\eps/4}\setminus S$, by the choice of
$e$ in (\ref{bigsmall}). Therefore, if $t\in [t_i,t_{i+1})$, and
if
\begin{equation}
\label{stir} \sup(e)\le \kappa\underline{\mathbf d}(\pi)
\end{equation}
then the choice of $\kappa$ gives
\begin{eqnarray}
\label{dgg}
 V(x_{i+1})-V(x_i)&= & V(x_{i+1})-V(\tilde x_\pi(t^-_{i+1}))+
 V(\tilde x_\pi(t^-_{i+1}))-V(\tilde x_i)\nonumber\\&+&
V(\tilde x_i)-V(x_i)\nonumber\\{}&\le & {\cal
L}_\eps|x_{i+1}-\tilde
x_\pi(t^-_{i+1})|-\frac{t_{i+1}-t_i}{4}V(\tilde
x_i)\nonumber\\
&+&{\cal L}_\eps|e(t_i)|\; \; ({\rm by}\; (\ref{circc}))\\&\le&
{\cal L}_\eps |e(t_i)|e^{L\delta} -\frac{t_{i+1}-t_i}{4}V(\tilde
x_i) +{\cal L}_\eps|e(t_i)|\; \; ({\rm by}\; (\ref{circstar}))
\nonumber\\{}& \le &
\frac{\lambda_-}{16}(t_{i+1}-t_i)-\frac{t_{i+1}-t_i}{4}V(\tilde
x_i)\; \; ({\rm by}\; (\ref{stir})) \nonumber\\{}&\le&
-\frac{t_{i+1}-t_i}{8}V(\tilde x_i)\; \; ({\rm by}\;
(\ref{check})) \nonumber\\{}&\le &
-\frac{t_{i+1}-t_i}{8}V(x_i)+\frac{t_{i+1}-t_i}{8}|e(t_i)| {\cal
L}_\eps\nonumber\\{}&\le &
-\frac{t_{i+1}-t_i}{8}V(x_i)+\frac{(t_{i+1}-t_i)^2}{16}\lambda_-\nonumber\\{}&\le
&
 -\frac{t_{i+1}-t_i}{16}V(x_i)\nonumber
 \end{eqnarray}
where we use \[t_{i+1}-t_i\; \le \; \overline{\mathbf d}(\pi)\;
\le\;  \delta<1\] to get the last inequality.  Set
\[J(t)=\frac{16}{16+t}\] for all $t\ge 0$. One can easily check that
$Q^\eps$ contains  the set \[S_V:=\{p: V(p)\le \max\{V(q): |q|\le
M+N\}\}.\]  In fact, $p\in S_V$ implies
\begin{eqnarray*}
|p|&\le & \overline{\alpha}\left(\max\{V(q): |q|\le M+N\}\right)\\
&= &
\max\{\overline\alpha\circ\underline\alpha^{-1}\circ\underline\alpha(V(q)):
|q|\le M+N\}\\
&\le &\overline\alpha\circ\underline\alpha^{-1}(M+N).
\end{eqnarray*}
In particular, $x_o\in S_V$. It follows from (\ref{dgg}) that if
none of $x_o, x_1, \ldots, x_j$ lie in $S^{\tilde\eps/16}$, then
\begin{eqnarray*}
V(x_1)-V(x_0) & \le &
-\frac{t_1}{16}V(x_j)\\
V(x_2)-V(x_1) &\le &
-\frac{t_2-t_1}{16}V(x_j)\\
& \vdots & \\
V(x_j)-V(x_{j-1}) &\le & -\frac{t_j-t_{j-1}}{16}V(x_j).
\end{eqnarray*}
Summing the preceding inequalities would then give
\[
V(x_j)-V(x_o) \; \le \;  -\frac{t_j}{16}V(x_j), \; \; \; \; \;
{\rm so}\; \;  \; \; \; V(x_j)\; \le\;  J(t_j)V(x_o)\; \; .
\]
Hence,
\[
V(x_i)\le J(t_i)V(x_o)\; \; \; {\rm for}\;  \; \; i=0,1,\ldots, j.
\]
By the choice of $\delta$ in (\ref{circlecross}), it would then
follow from (\ref{shrp}) that
\[
V(x_\pi(t))\le J(t) V(x_o)+\frac{\tilde\eps}{8}
\]
up to the least time $t$ at which $x_\pi(t)\in S^{\tilde\eps/16}$.
Hence, for such $t$, the choice of $\tilde\eps$ (see
(\ref{stars})) gives
\[
\begin{array}{ccl}
|x_\pi(t)| &\le&
\overline{\alpha}\left(J(t)V(x_o)+\frac{\tilde\eps}{8}\right)\\
&\le & \overline{\alpha}\left(J(t)V(x_o)\right)+\frac{\eps}{8}
\end{array}.
\]
 On the other
hand, (\ref{dgg}) also shows that if $x_\pi(t)\in
S^{\tilde\eps/8}$ for some $t$, then
\begin{equation}
\label{neww} |x_\pi(s)|\le \overline{\alpha}\circ
\underline\alpha^{-1}(N)+\eps \; \; \forall s\ge t. \end{equation}
Indeed, let $s_1$ be the first sample time above such a time $t$.
Assume $x_\pi(t)\not\in \eps {\cal B}_n$.  By (\ref{shrp}),
$x_\pi(s_1)\in S^{\tilde\eps/4}$ and $x_\pi(s_1)\not\in
\frac{\eps}{2}{\cal B}_n$.  Therefore, there exists $p\in S$ for
which
\begin{eqnarray}
V(x_\pi(s_1)) & = & V(x_\pi(s_1))-V(p)+V(p)\nonumber\\
&\le & {\cal L}_\eps
\frac{\tilde\eps}{4}+\underline{\alpha}^{-1}(N)\nonumber.
\end{eqnarray}
In fact, we can pick $p=x_\pi(s_1)$ if $x_\pi(s_1)\in S$ and $p\in
\partial S$ otherwise, so $p\not\in \frac{\eps}{2}{\cal B}_n$.
It follows from (\ref{shrp}) and (\ref{dgg}) that for the next
sample time $s_i$, we either have  $x_\pi(s_i)\in S^{\tilde
\eps/8}$, or else we have
\[V(x_\pi(s_i))\le {\cal
L}_\eps\frac{\tilde\eps}{4}+\underline{\alpha}^{-1}(N).\] In the
first case, \[|x_\pi(s_i)|\le \overline{\alpha}\circ
\underline{\alpha}^{-1}(N)+\frac{\eps}{8},\] while in the second
case,
\[|x_\pi(s_i)|\le \overline{\alpha}\left(\frac{\tilde\eps{\cal L}_\eps}{4}
+\underline{\alpha}^{-1}(N)\right)\le
\overline{\alpha}\circ\underline{\alpha}^{-1}(N)+\frac{\eps}{8},\]
by the choice of $\tilde\eps$.  If $x_\pi(s_i)\not\in S^{\tilde
\eps/16}$, then $V(x_\pi(s_{i+1}))\le V(x_\pi(s_{i}))$ (by
(\ref{dgg})), so the preceding argument also gives
\[|x_\pi(s_{i+1})|\le \overline{\alpha}\circ
\underline{\alpha}^{-1}(N)+\frac{\eps}{8}.\] By repeating this
argument for subsequent sample times, the assertion (\ref{neww})
then follows from (\ref{shrp}). Defining $\beta\in {\cal KL}$ and
$\gamma\in {\cal K}_\infty$ by
\begin{equation}
\label{klk}
\beta(s,t)=\overline{\alpha}\left(\underline{\alpha}^{-1}(s)J(t)\right),
\; \; \gamma(s)=\overline{\alpha}\circ\underline\alpha^{-1}(s),
\end{equation}
it follows that (\ref{goal}) holds for all $x_o\in M\bar {\cal
B}_n$, $u\in {\cal M}^m_N$, $\pi\in {\rm Par}(\delta)$, and  $e\in
{\cal O}$ for which $\sup(e)\le \kappa \underline{\mathbf
d}(\pi)$. Therefore, Theorem~\ref{thm2} will follow once we  check
(\ref{circstar}), which is a consequence of Gronwall's Inequality.

To this end, notice that if $x_i\in Q^{\eps/4}$, then
\[
|x_\pi(t)-\tilde x_\pi(t)|\le |x_i-\tilde x_i|+ \int_{t_i}^{t}
({\cal L}_f|x_\pi(s)-\tilde x_\pi(s)|+R{\cal L}_G|x_\pi(s)-\tilde
x_\pi(s)|)\, {\rm d}s
\]
for all $t\in[t_i,t_{i+1})$, where we are using the constants in
(\ref{constants}).   It follows from Gronwall's Inequality that
\[
|x_\pi(t)-\tilde x_\pi(t)|\le |x_i-\tilde x_i|
e^{L|t_i-t_{i+1}|}\le |x_i-\tilde x_i| e^{L\overline{\mathbf
d}(\pi)}\le |e(t_i)|e^{L\delta}
\]
for all $t\in [t_i,t_{i+1})$, which is (\ref{circstar}). This
proves Theorem~\ref{thm2}.

We turn next to Theorem~\ref{thm1}. We need to show the ISS
property  (\ref{iss}) for all Euler solutions $x(t)$ of
(\ref{store}).  We will actually prove the slightly stronger
version of the theorem for generalized Euler solutions, as
asserted in Remark~\ref{euleremark}. To this end, choose $u\in
{\cal M}^m$, $x_o\in \R^n$, and $\eps>0$.  Using our previous
conclusion that (\ref{two}) is ISS for sampling solutions, we can
let
\[\delta_\eps=\delta\left(\eps, |x_o|, |u|_\infty\right)\; \;
\; {\rm and}\; \; \;  \kappa_\eps=\kappa\left(\eps, |x_o|,
|u|_\infty\right)\] be the constants from Definition~\ref{Defn3}.
Let $x(t)$ be a generalized Euler solution of (\ref{store}), and
let $\pi_r$, $u_r$, and $e_r$ satisfy the requirements of the
generalized Euler solution definition. It follows from the
definition that there is an $\bar r\in \N$  such that
\[\overline{\mathbf d}(\pi_{r})\le \delta_\eps,\; \; \;
\sup(e_{r})\le \kappa_\eps\underline{\mathbf d}(\pi_{r})\] for all
$r\ge \bar r$. It then follows from (\ref{goal}) that
\begin{equation}
\label{goali} |x_{\pi_{r}}(t;x_o, u_{r}, e_{r})|\le
\beta(|x_o|,t)+\gamma(|u|_\infty)+\eps
\end{equation}
for all $t\ge 0$ and $r\ge \bar r$, where $\beta$ and $\gamma$ are
in (\ref{klk}).  The ISS condition (\ref{iss}) now follows by
passing to the limit in (\ref{goali}) as $r\to \infty$, since
$\eps>0$ was arbitrary. This concludes the proof of Theorem
\ref{thm1}.

\section{Stabilization of the Nonholonomic Integrator}
\label{fbrk} In this section, we illustrate how the feedback
constructed in $\S$\ref{proofs} can be used to stabilize
Brockett's nonholonomic integrator control system (see \cite{B83,
LR03, S99}).  We will also use the nonholonomic integrator to
compare our feedback construction to the feedbacks from
\cite{S89a, S90}. The nonholonomic integrator was introduced in
\cite{B83}, as an example of a system that cannot be stabilized
using continuous feedback.  It is well-known that if the state
space of a system contains obstacles (e.g., if the state space is
$\R^2\setminus(-1,1)^2$, and therefore has a topological obstacle
around the origin), then it is impossible to stabilize the system
using continuous feedback.  In fact, this is a special case of a
theorem of Milnor, which asserts that the domain of attraction of
an asymptotically stable vector field must be diffeomorphic to
Euclidean space, and therefore cannot be the complement
$\R^2\setminus(-1,1)^2$ (see \cite{S98}).

Brockett's example illustrates how, even if we assume that the
state evolves in Euclidean space, similar obstructions to
stabilization may occur.  These obstructions are not due to the
topology of the state space, but instead arise from ``virtual
obstacles'' that are implicit in the form of the control system
(see \cite{S99}).  Such obstacles occur when it is impossible to
move {\em instantly} in some directions, even though it is
possible to move {\em eventually} in every direction
(``nonholonomy'').  This gives rise to Brockett's criterion (see
\cite{B83}), which is a necessary condition for the existence of a
continuous stabilizer, in terms of the vector fields that define
the system (see \cite{S98, S99, St02}).  The nonholonomic
integrator does not satisfy Brockett's criterion, and therefore
cannot be stabilized by continuous feedbacks.

The physical model for Brockett's example is as follows. Consider
a three-wheeled shopping cart whose front wheel acts as a castor.
The state variable is $(x_1,x_2,\theta)^T$, where $(x_1,x_2)^T$ is
the midpoint of the rear axle of the cart, and $\theta$ is the
cart's orientation.  The front wheel is free to rotate, but there
is a ``non-slipping'' constraint that $(\dot x_1, \dot x_2)^T$
must always be parallel to $(\cos(\theta),\sin(\theta))^T$.  This
gives the equations
\begin{equation}
\label{slip}
\begin{array}{lll}
\dot x_1 & = & v_1\cos(\theta)\\
\dot x_2 & = & v_1\sin(\theta)\\
\dot \theta&=&v_2
\end{array},
\end{equation}
where $v_1$ is a ``drive'' command and $v_2$ is a steering
command. Using the feedback transformation
\[\begin{array}{l}
z_1:=\theta,\; z_2:=x_1\cos(\theta)+x_2\sin(\theta),\;
z_3:=x_1\sin(\theta)-x_2\cos(\theta)\\  u_1:=v_2,\;
u_2:=v_1-v_2z_3
\end{array}
\] followed by a second
transformation, brings the equations (\ref{slip}) into the
 form
\begin{equation}
\label{broc}
\begin{array}{lll}
\dot x_1&=&u_1\\
\dot x_2&=&u_2\\
 \dot x_3&=&x_1u_2-x_2u_1
 \end{array},
 \end{equation}
which is called the nonholonomic integrator control system.

One can show (see  \cite{LS99}) that (\ref{broc}) is a GAC system.
However, since Brockett's  condition is not satisfied for
(\ref{broc}), the system has no continuous stabilizer. While there
does not exist a $C^1$ CLF for the system (\ref{broc}) (see
\cite{LS99}), it is now well-known that every GAC system admits a
continuous CLF (see  \cite{S83}). In fact, it was shown in
\cite{LR03} that the nonholonomic integrator dynamics (\ref{broc})
has the  nonsmooth CLF
\begin{equation}
\label{newclf} V(x)=\max\left\{\sqrt{x^2_1+x^2_2},
|x_3|-\sqrt{x^2_1+x^2_2}\right\}
\end{equation}
which is semiconcave outside the cone $x^2_3=4(x^2_1+x^2_2)$
 (see \cite{R02} for a
detailed discussion of some special properties of this CLF). For
the special case of the dynamics (\ref{broc}) and CLF
(\ref{newclf}), the feedback $K=K_1+K_2$ we constructed in
$\S$\ref{proofs} is as follows.

To simplify notation, we use the radius
$r(x):=\sqrt{x^2_1+x^2_2}$.  We also use    the sets
\[\begin{array}{lcl}
S_o&=&\{x\in\R^3:  \; x_3\ne 0, \; r(x)=0\}\\
S_+&=&\{x\in\R^3: \; x^2_3\ge 4r^2(x)>0\}\\
S_-&=&\{x\in\R^3: \; x^2_3< 4r^2(x)\}
\end{array}
\]which form a partition of $\R^3\setminus\{0\}$.  Notice that
$V(x)=r(x)$ on $S_-$,  and also that $V(x)=|x_3|-r(x)$ on
$\R^3\setminus S_-$. To find our selection $\zeta(x)\in
\partial_LV(x)$, we first choose  $\zeta(0)=0$,
and $\zeta(x)=(0,-1,\sgn\{x_3\})^T$ for all  $x\in S_o$. Using the
notation of (\ref{newww}), this gives
\begin{equation}
\label{rou} b(x)=\left\{
\begin{array}{ll}
\left(\,  \ds -x_2\sgn\{x_3\}-x_1/r(x), \; \ds
x_1\sgn\{x_3\}-x_2/r(x) \, \right)^T, &
x\in S_+\\
\ds\left(\, x_1/r(x), \; x_2/r(x)\, \right)^T, & x\in S_-\\
\left( 0, -1\right)^T, & x\in S_o
\end{array}\right.
\end{equation}
and $b(0)=0$.  Notice that $1\le |b(x)|^2\le r^2(x)+1$ for all
$x\ne 0$.  We also have
\[
K_1(x)=\left\{
\begin{array}{ll}
\! \! \mu_1(x) \left(\, -x_2\sgn\{x_3\}-x_1/r(x), \;
x_1\sgn\{x_3\}-x_2/r(x) \, \right)^T, &
\! \! x\in S_+\\
\ds -\left(x_1, x_2\right)^T, & \! \! x\in S_-\\
\left( 0, |x_3|\right)^T, & \! \! x\in S_o
\end{array}\right.
\]
with $K_1(0)=0$, where we have set
\[\mu_1(x):=
\frac{r(x)-|x_3|}{r^2(x)+1}.
\]
 In this case, we have taken
\[K_1(x)=-b(x)V(x)/|b(x)|^2\] for $x\ne 0$, where $b(x)$ is defined
in (\ref{rou}), and $K_1$ is continuous at the origin.
 On
the other hand, our feedback $K_2$ from (\ref{newww}) becomes
\[
K_2(x)=-\left\{
\begin{array}{ll}
 \left(\, \mu_2(x_1,-x_2,x),\;
\mu_2(x_2,x_1,x)\, \right)^T, & x\in S_+\\
r(x)\left(\, \sgn\{x_1\},\,
\sgn\{x_2\}\, \right)^T, & x\in S_-\\
|x_3|\left( 0,\;  -1 \right)^T,& x\in S_0
\end{array}\right.
\]
with $K_2(0)=0$, where we have set
\[
\mu_2(a,b,x)\; :=\; (|x_3|-r(x))\, \sgn\{\, b\, r(x)\,
\sgn\{x_3\}-a\, \}.
\]
Since $V$ is semiconcave on $\Omega:=\R^3\setminus {\rm bd}(S_-)$,
the argument from $\S$\ref{proofs} applies to sampling solutions
that satisfy the additional requirement that $\tilde x_\pi(s)\in
\Omega$ for all $s\ge 0$.
 It follows from the proof of Theorem~\ref{thm2}
that the nonholonomic integrator system (\ref{broc}) can be
stabilized for both actuator errors and small observation errors
(for this restricted set of sampling solutions), using the
combined feedback $K=K_1+K_2$.

\begin{remark}\rm
 In this example, we chose to work with
the  CLF (\ref{newclf}) because it has been explicitly proven in
\cite{LR03} to be a CLF for the control system (\ref{broc}). The
example illustrates how to extend our results to more general
CLF's that may not be semiconcave on $\R^3\setminus\{0\}$. For
such cases, the ISS estimates hold for those sampling solutions
that remain in the domain of semiconcavity of the CLF.  On the
other hand, we let the reader prove that the nonholonomic
integrator system also has the CLF \[\tilde V(x) =
\left(\sqrt{x_1^2+x_2^2} - |x_3|\right)^2 + x_3^2,\] which is
semiconcave on $\R^3\setminus\{0\}$ (as the sum of the smooth
function $x_1^2+x^2_2+2x^2_3$ and a semiconcave function).
Therefore, if we use $\tilde V$ to form our feedbacks, instead of
the CLF (\ref{newclf}), then our theorems apply directly, without
any state restrictions on the sampling solutions.
\end{remark}

\begin{remark} \label{origin}\rm  The results in \cite{S89a} designed
feedbacks that make $C^o$-stabilizable systems ISS with respect to
actuator errors. For the case of $C^o$-stabilizable systems, a
smooth (i.e., $C^\infty$) Lyapunov function is known to exist (see
\cite{A83}). In \cite{S89a}, the system was rendered ISS using the
feedback
\begin{equation}
\label{trick} \hat K(x):=-L_GV(x)=-\nabla V(x)G(x),
\end{equation}
where $V$ is a smooth CLF for the dynamics (\ref{dyn}). In that
case, (\ref{trick}) is continuous at the origin. However, in the
more general situation where the system is merely GAC, there may
not exist a smooth Lyapunov function, so $V$ must be taken to be
nonsmooth.  In this case, the use of the nonsmooth analogue
\begin{equation}
\label{analogue} \tilde K(x):=-\zeta(x)G(x)\end{equation} of
(\ref{trick}) (where $\zeta(x)\in
\partial_LV(x)$ for all $x\ne 0$) could give rise to a feedback
that would not be continuous at the origin.
 For example, if we use the nonholonomic integrator (\ref{broc})
 and the CLF (\ref{newclf}), then $\tilde K$ takes the values
\[
\tilde K\left((\eps,\eps,0)^T\right)=-\left(\frac{1}{\sqrt{2}},
\frac{1}{\sqrt{2}}\right)^T,\; \; \tilde K\left((\eps,\eps,
3\sqrt{2}\eps)^T\right)= \left( \frac{1}{\sqrt{2}},
\frac{1}{\sqrt{2}}\right)^T+\eps(1,-1)^T
\]
so $\tilde K$ is discontinuous at the origin.   On the other hand,
our choice of $K_2$ is automatically continuous at the origin. \er
\br{scp} \rm Under the additional hypothesis that (\ref{dyn})
satisfies the small control property (see  \cite{S98}), the system
can be stabilized by a feedback that is  continuous at the origin
(see \cite{R02}). More precisely, suppose there exists a
semiconcave CLF $V$ satisfying the following:  For each $\eps>0$,
there exists $\delta=\delta(\eps)>0$ such that $0<|x|\le \delta$
implies
\[
\exists u_x\in \eps {\cal B}_m\; \; {\rm s.t.}\; \; \forall
\zeta\in
\partial_PV(x),\; \; \langle\zeta, f(x)+G(x)u_x\rangle\le -V(x).
\]
Then the system can be rendered globally asymptotically stable
(GAS) by a feedback that is  continuous at the origin  (see
\cite{R02}).  For the case of the nonholonomic integrator
(\ref{broc}), the system is GAS under our feedback $K_1$, which is
continuous at the origin, so our total feedback $K=K_1+K_2$ is
continuous at the origin as well.\er

\section{ISS for Fully Nonlinear GAC Systems}
\label{rks} We conclude with an extension of our results for fully
nonlinear GAC systems
\begin{equation}\label{ffeq}
  \dot x=f(x,u)
\end{equation}
where we assume for simplicity that the observation error $e$ in
the controller is zero. We assume throughout this section that
\[f:\R^n\times \R^m\to \R^n: (x,u)\mapsto f(x,u)\] is continuous and
locally Lipschitz in $x$ uniformly on compact subsets of
$\R^n\times \R^m$ and $f(0,0)=0$. It is natural to ask whether
these hypotheses are sufficient for the existence of a continuous
feedback $K(x)$ for which
\begin{equation}\label{ship}
  \dot x=f(x,K(x)+u)
\end{equation}
is ISS for Euler solutions.  However, one can easily construct
examples  for which  such feedbacks cannot exist.  Here is an
example from \cite{S90} where this situation occurs.  Consider the
GAC system $\dot x=-x+u^2x^2$  on $\R$.  If $K(x)$ is any
continuous feedback for which
\begin{equation}
\label{exx} \dot x=-x+(K(x)+u)^2x^2
\end{equation}
is ISS, then  $|K(x)|<x^{-1/2}$ for sufficiently large $x>0$. It
follows that the solution of \[\dot x=-x+(K(x)+1)^2x^2\] starting
at $x(0)=4$ is unbounded.  Therefore, there does not exist a
continuous feedback $K$ for which (\ref{exx}) is ISS. On the other
hand, one {\em can} find a (possibly discontinuous) feedback that
makes (\ref{ffeq}) ISS. We use the following weaker sense of ISS
for fully nonlinear systems that was introduced in \cite{S90}:

\bd{Defn6} We say that (\ref{ffeq}) is {\rm  input to state
stabilizable in the weak sense} provided there exist a feedback
$K$, and an $m\times m$ matrix $G$ of continuously differentiable
functions which is invertible at each point, such that
\[\dot x=F(x,K(x),u)\] is ISS for sampling and
Euler solutions, where $F(x,p,u)=f(x,p+G(x)u)$. \ed

We will prove the following: \bp{var} If (\ref{ffeq}) is GAC, then
(\ref{ffeq}) is also  input to state stabilizable in the weak
sense.\ep

\begin{proof}
 We  modify the proof from $\S$\ref{proofs}.  We
define $V$, $\zeta$, $\underline{\alpha}$, $\overline{\alpha}$,
and $K_1$ as in the proof of Theorem~\ref{thm2}, except  we use
the fully nonlinear dynamics $h=f$ from (\ref{ffeq}). Next we
follow the proof of the main result in \cite{S90}, with the
following modifications. Define the (possibly discontinuous)
function ${\cal D}$ by
\begin{equation}\label{dee}
{\cal D}(s,r)=\sup\left\{\langle\zeta(x),f(x,K_1(x)+p)\rangle
+\frac{V(x)}{2}: |x|=s, |p|=r \right\}.
\end{equation}
For any interval $I$ of the form $[i,i+1]$, or of the form
$[\frac{1}{i+1},\frac{1}{i}]$, for $i\in\N$, one can find
$r=r(i)>0$ such that $s\in I$ implies ${\cal D}(s,b)<0$ for all
$b\in [0,r]$.  This follows from the positive definiteness of $V$,
the local Lipschitzness of $f$,  and the local boundedness of
$\partial_PV$ on compact subsets of $\R^n\setminus\{0\}$.

The argument of \cite{S90} therefore gives $\alpha_4\in {\cal
K}_\infty$ and a smooth, everywhere invertible matrix-valued
function $G:\R^n\to\R^{m\times m}$ satisfying the following:  If
\begin{equation}\label{ifn} |x|> \alpha_4(|u(s)|_\infty),\end{equation} then for a.e. $t\ge 0$,
\[
\left\langle \zeta(x), f(x, K_1(x)+G(x)u(t))\right\rangle
+\frac{V(x)}{2}\; \le\; {\cal D} (|x|, |G(x)u(s)|_\infty)<0.
\]
(See Remark~\ref{expe} for a characterization of the set of
matrices $G$ for which ISS can be expected, in terms of ${\cal
D}$.) We can evidently assume that $\alpha_4(s)\ge s$ for all
$s\ge 0$ (e.g., by replacing $\alpha_4(s)$ by
$\max\{\alpha_4(s),s\}$, which makes the condition (\ref{ifn})
more restrictive). Fix $M$, $N$, $\eps\in (0,M)$, $u\in {\cal
M}^m_N$ and $x(t)=x_\pi(t)$ as before, with $e=0$. Define the
compact sets
\[S:=\left\{x\in\R^n: V(x)\le \underline{\alpha}^{-1}\circ
\alpha_4(N)\right\},\; \; \; Q=\left\{(\overline{\alpha}\circ
\underline{\alpha}^{-1}(M+\alpha_4(N))+1)\bar {\cal
B}_n\right\}\setminus \eps {\cal B}_n.\]  Notice that $S\subseteq
Q^\eps$.  We choose $\tilde\eps$ as before, and we choose
$\delta=\delta(\eps,M,N)$, satisfying (\ref{circlecross}), such
that if $\overline{\mathbf d}(\pi)<\delta$, then
\begin{equation}
\label{shrp2}
 |x_\pi(t)-x_i|\; \le\;
 \min\left\{\mu,\frac{\tilde \eps}{16(1+{\cal L}_\eps)},
\sqrt{\frac{\lambda_-}{8\sigma}(t-t_i)}\right\}
\end{equation}
for all indices $i$ such that $x_i\in Q^\eps$ and all
$t\in[t_i,t_{i+1}]$, where $\sigma$ and  $\mu$ are as defined
before, and $\lambda_-=\min\{V(x): x\in Q^{\eps/4}\}$. Reducing
$\delta$ as necessary, we can assume
\[
\ds \Vert\zeta(x_i)\cdot [
f(x_i,K_1(x_i)+G(x_i)u(s))-f(x(s),K_1(x_i)+G(x(s))u(s))]\Vert_{[t_i,t_{i+1}]}
\le \frac{\lambda_-}{8}
\]
 for all indices $i$ satisfying
 $x_i\in Q^{\eps/2}$. Reasoning as in the earlier proof gives
\[
V(x_\pi(t))-V(x_i)\le -(t-t_i)\frac{V(x_i)}{16}\; \; \forall t\in
[t_i, t_{i+1}]
\]
for all $i$ such that $x_i\in Q^{\eps/4}\setminus S$.  The
remainder of the proof is as before, except with
$\overline{\alpha}\circ \underline{\alpha}^{-1} (N)$ replaced by
$\overline{\alpha}\circ \underline{\alpha}^{-1} (\alpha_4(N))$,
and with $\overline{\alpha}\circ \underline{\alpha}^{-1} (s)$
replaced by $\overline{\alpha}\circ \underline{\alpha}^{-1}
(\alpha_4(s))$ in the definition of $\gamma$.  This proves
Proposition~\ref{var}. \qquad\end{proof}

\begin{remark}\label{expe}\rm The statement of
Proposition~\ref{var} is an existence
result in terms of the invertible matrix $G$.  However, we can
strengthen the proposition by using the function ${\cal D}$ in
(\ref{dee}) to  characterize the class of $G$  for which ISS can
be expected, as follows. Following \cite{S90}, we first choose
strictly decreasing sequences $\{r_i\}$ and $\{r'_i\}$ of positive
numbers such that $0<r_{i+1}<r'_i<r_i$ for all $i\in\N$, and such
that
\[
{\cal D}(s,r)<0 \; \; \; \forall (s,r)\in \left([i,i+1] \times
[0,r_i]\right)\cup  \left([1/(i+1),1/i]\times [0,r'_i] \right)
\]
for all $i\in \N$.  The existence of these sequences follows from
the argument we gave in the proof of the proposition. Define
 $\rho:[0,\infty)\to[0,\infty)$ by setting:
\begin{itemize}\item[($\rho1$)]$\rho(s)=r_k$ for all $s\in
[k,k+1)$ and $k\in \N$; \item[($\rho$2)]$\rho(s)=r'_k$ for all
$s\in [1/(k+1),1/k)$ and $k\in \N$; and
\item[($\rho$3)]$\rho(0)=0$.\end{itemize} We then choose any
smooth function $g: [0,\infty)\to (0,\infty)$ satisfying:
\begin{itemize}\item[($g$1)]$g(s)=1$ for all $s\in [0,1]$;
\item[($g$2)] $g(s)\le \rho(s)/s$ for all $s\ge 2$; and
\item[($g$3)]$g(s)\le 1$ for all $s\ge 0$.\end{itemize} The
existence of such a function $g$ follows from exactly the same
argument used in \cite{S90}.  It then also follows from the
argument of \cite{S90} that we can satisfy the conditions of the
proposition by choosing $G(\xi)=g(|\xi|)I$. \er

Proposition~\ref{var} allows us to characterize GAC for fully
nonlinear systems in terms of feedback equivalence, as follows.
First recall that two systems $\dot x=f(x,u)$ and $\dot x=h(x,u)$,
evolving on $\R^n\times \R^m$, are called {\it feedback
equivalent} provided there exist a locally bounded function
$K:\R^n\to \R^m$ and an everywhere invertible function
$G:\R^n\to\R^{m\times m}$ for which \[h(x,u)=f(x,K(x)+G(x)u)\] for
all $x\in \R^n$ and $u\in\R^m$;  in this case, we also say $\dot
x=f(x,u)$ is feedback equivalent to (\ref{feed}) with $e\equiv 0$
and $F(x,p,u)=f(x,p+G(x)u)$. The following elegant statement
follows directly from Proposition~\ref{var}:

\bc{corlast} The fully nonlinear system (\ref{ffeq}) is GAC if,
and only if, it is feedback equivalent to a system which is ISS
for sampling and Euler solutions. \ec

\begin{remark}\rm Although, as shown by the counterexample (\ref{exx}), it
is in general impossible to obtain input to state stabilization
(in the non-weak sense) for systems that are not affine in
controls, it is still the case that for some restricted classes of
systems this objective can be attained, under appropriate
neutral-stability assumptions on the dynamics. One such class is
that of systems in which the input appears inside a saturation
nonlinearity, such as $\dot x = f(x,u) = f_0(x) + g(x)\sigma(u)$.
The papers~\cite{liu-chitour-sontag1} and
\cite{liu-chitour-sontag2} (see~\cite{sussmann-sontag-yang} for an
application of these results to the recursive design of
stabilizers for a large class of systems with saturation) as well
as~\cite{chitour} and~\cite{lin} dealt with such questions, for
systems that are linear in the absence of the saturation (the
$f_0$ and $g$ vector fields are linear and constant,
respectively), while~\cite{bao-lin} obtained similar results for
more general nonlinear systems. \er

\bibliographystyle{siam}

\end{document}